\documentclass[12pt,reqno]{amsart}

\usepackage{amscd}

\theoremstyle{plain}
\newtheorem{lemma}{Lemma}[section]
\newtheorem{theorem}[lemma]{Theorem}
\newtheorem{proposition}[lemma]{Proposition}
\newtheorem{corollary}[lemma]{Corollary}
\newtheorem{conjecture}[lemma]{Conjecture}

\theoremstyle{definition}
\newtheorem{definition}[lemma]{Definition}

\newtheorem{remark}[lemma]{Remark}
\newtheorem{remarks}[lemma]{Remarks}
\newtheorem{example}[lemma]{Example}

\newcommand{\N}{{\mathbb N}}

\newcommand{\Q}{{\mathbb Q}}
\newcommand{\C}{{\mathbb C}}

\renewcommand{\P}{{\mathbb P}}

\newcommand{\MM}{{\mathcal M}}
\newcommand{\OO}{{\mathcal O}}
\newcommand{\Nu}{{\mathcal V}}
\newcommand{\KK}{{\mathcal K}}

\newcommand{\TT}{{\mathcal T}}

\newcommand{\aaa}{{\bf a}}

\newcommand{\mmm}{{\bf m}}

\newcommand{\www}{\widetilde}
\newcommand{\tm}{{\mathcal T}_M}
\newcommand{\paa}{\partial}
\newcommand{\oooo}{\overline}

\newcommand{\nnn}{\nabla}

\DeclareMathOperator{\Ann}{Ann}

\DeclareMathOperator{\Lie}{Lie}

\DeclareMathOperator{\id}{id}

\DeclareMathOperator{\trace}{trace}

\DeclareMathOperator{\Gr}{Gr}

\DeclareMathOperator{\Specan}{Specan}

\begin{document}

\title{Variance of the spectral numbers}

\author{
Claus Hertling}

\address{Claus Hertling\\
Mathematisches Institut der Universit\"at Bonn\\
Beringstra{\ss}e 1, 53115 Bonn, Germany}

\email{hertling\char64 math.uni-bonn.de}


\begin{abstract}
A formula for the variance of the spectrum of a quasihomogeneous 
singularity is proved, using the G-function of a semisimple 
Frobenius manifold.
\end{abstract}

\maketitle









\section{Introduction}
\setcounter{equation}{0}

The spectrum of an isolated hypersurface singularity
$f:(\C^{n+1},0)\to (\C,0)$ is an important discrete invariant of the 
singularity. Its main properties have been established by Steenbrink
and Varchenko. It consists of $\mu$ rational numbers
$\alpha_1,...,\alpha_\mu$ with 
$-1<\alpha_1\leq ... \leq \alpha_\mu < n$ and 
$\alpha_i+\alpha_{\mu+1-i}=n-1$.
The numbers $e^{-2\pi i\alpha_1},...,e^{-2\pi i\alpha_\mu}$ are the
eigenvalues of the monodromy. 

The spectral numbers come from a Hodge filtration on the cohomology of
a Milnor fiber \cite{St1} (cf. chapter 4) or, more instructively, 
from the Gau{\ss}-Manin connection of $f$ \cite{AGV}.
They satisfy a semicontinuity property for deformations of $f$ and are 
related to the signature of the intersection form.

In the case of a quasihomogeneous singularity $f$ of weighted degree 1 
with weights $w_0,...,w_n\in (0,\frac{1}{2}]\cap \Q$
they can be calculated easily \cite{St2}\cite{AGV}:
Then the Jacobi algebra 
$\OO_{\C^{n+1},0}/(\frac{\paa f }{\paa x_0},...,\frac{\paa f}{\paa x_n})
=:\OO/J_f$
has a natural grading $\OO/J_f = \bigoplus_\alpha (\OO/J_f)_\alpha$, and
\begin{eqnarray}
\sharp (i\ |\ \alpha_i=\alpha) &=& \dim (\OO/J_f)_{\alpha-\alpha_1}\ ,\\
\alpha_1 &=& -1+\sum_{i=0}^n w_i\ .
\end{eqnarray}

The main result of this paper is a new formula concerning the distribution
of the spectral numbers. Because of the symmetry 
$\alpha_i+\alpha_{\mu+1-i}=n-1$ one can consider
$\frac{n-1}{2}$ as their {\it expectation value}. Then
$\frac{1}{\mu} \sum_{i=1}^\mu \left(\alpha_i-\frac{n-1}{2}\right)^2$
is their {\it variance}.

\begin{theorem}
In the case of a quasihomogeneous singularity the variance is
\begin{eqnarray}
\frac{1}{\mu} \sum_{i=1}^\mu \left(\alpha_i-\frac{n-1}{2}\right)^2
= \frac{\alpha_\mu-\alpha_1}{12}\ .
\end{eqnarray}
\end{theorem}

\begin{conjecture}
For any isolated hypersurface singularity
\begin{eqnarray}
\frac{1}{\mu} \sum_{i=1}^\mu \left(\alpha_i-\frac{n-1}{2}\right)^2
\leq \frac{\alpha_\mu-\alpha_1}{12}\ .
\end{eqnarray}
\end{conjecture}

The proof of (1.3) uses two things: 

1) The deep result of K. Saito 
\cite{SK1}\cite{SK3} and M. Saito \cite{SM} that the base space of a
semiuniversal unfolding of an isolated hypersurface singularity $f$
can be equipped with the structure of a Frobenius manifold with
discrete invariants related to the spectrum of $f$.

2) The {\it G-function} of a semisimple Frobenius manifold, which was
defined by Dubrovin and Zhang \cite{DZ1} and independently by
Givental \cite{Gi}.

A more elementary and much broader version of the construction of K. Saito
and M. Saito has been given in \cite{He4} (chapter 6). In chapter 4 we will
state the result more precisely. The definition of a Frobenius manifold
is given in chapter 3.

The G-function of a semisimple Frobenius manifold is a fascinating function
with several origins. One is the {\it $\tau$-function} of the isomonodromic
deformations, which are associated to such a Frobenius manifold.
Another is, that in the case of a semisimple Frobenius manifold coming 
from quantum cohomology the G-function is the genus one Gromov-Witten 
potential. More remarks on this, the definition, and some properties of the
G-function are given in chapter 6.

In order to establish the fact that the G-function extends holomorphically
over the {\it caustic} in the singularity case (Theorem 6.3),
we need the {\it socle field} of a Frobenius manifold (chapter 5)
and some facts on {\it F-manifolds} (chapter 2) from \cite{He3}. 
If one forgets the metric of a Frobenius manifold 
one is left with an F-manifold.

Theorem 1.1 is proved in chapter 7.

Conjecture 1.2 is based only on a few examples. I would appreciate a proof
as well as counterexamples, also an elementary proof of Theorem 1.1, 
and applications, for example on deformations of singularities or on their
topology.

\section{F-manifolds}
\setcounter{equation}{0}

First we fix some notations: 

1) In the whole paper $M$ is a complex manifold of dimension $m\geq 1$
(with $m=\mu$ in the singularity case) with holomophic tangent bundle
$TM$, sheaf $\tm$ of holomorphic vector fields, and sheaf $\OO_M$ of 
holomorphic functions.

2) A $(k,l)$-tensor is an $\OO_M$-linear map 
$\tm^{\otimes k}\to \tm^{\otimes l}$.
The Lie derivative $\Lie_X T$ of it by a vector field $X\in \tm$ is again
a $(k,l)$-tensor. For example, a vector field $Y\in \tm$ yields a 
$(0,1)$-tensor $\OO_M\to \tm$, $1\mapsto Y$,
with $\Lie_XY = [X,Y]$.

3) If $\nnn$ is a connection on $M$ then the covariant derivative 
$\nnn_X T$ of a $(k,l)$-tensor by a vector field $X$ is again a 
$(k,l)$-tensor. As $\nnn_X T$ is $\OO_M$-linear in $X$ (contrary to 
$\Lie_X T$) $\nnn T$ is a $(k+1,l)$-tensor.

4) A {\it multiplication} $\circ$ on the tangent bundle $TM$ of a manifold $M$
is a symmetric and associative $(2,1)$-tensor $\circ :\tm\otimes \tm\to \tm$.
It equips each tangent space $T_tM$, $t\in M$, with the structure of a 
commutative and associative $\C$-algebra. We will be interested only in 
a multiplication with a global unit field $e$.

5) A {\it metric} $g$ on a manifold $M$ is a symmetric and nondegenerate 
$(2,0)$-tensor. It equips each tangent space with a symmetric and
nondegenerate bilinear form. Its Levi-Civita connection $\nnn$ is the 
unique connection on $TM$ which is torsion free, i.e.
$\nnn_X Y-\nnn_Y X = [X,Y]$, and which satisfies $\nnn g=0$, i.e. 
$X\ g(Y,Z)= g(\nnn_X Y,Z)+g(Y,\nnn_X Z)$.

\bigskip
The notion of an F-manifold was defined first in \cite{HM} (cf. \cite{Man}
I \S 5). It is studied extensively in \cite{He3}.

\begin{definition}
a) An {\it F-manifold} $(M,\circ ,e)$ is a manifold $M$ together with a 
multiplication $\circ$ on the tangent bundle  and a global unit field $e$
such that the multiplication satisfies the following 
{\it integrability condition}:
\begin{eqnarray}
\forall \ X,Y\in \tm \qquad
\Lie_{X \circ Y}(\circ)=X \circ \Lie_Y(\circ) + Y \circ\Lie_X(\circ)\ .
\end{eqnarray}
b) Let $(M,\circ,e)$ be an F-manifold. An {\it Euler field of weight} 
$c\in \C$ is a vector field $E\in \tm$ with 
\begin{eqnarray}
\Lie_E(\circ) = c\cdot \circ\ .
\end{eqnarray}
An Euler field of weight 1 is simply called an {\it Euler field}.
\end{definition}

One reason why this is a natural and good notion is that any Frobenius 
manifold is an F-manifold (\cite{HM}, \cite{He3} chapter 5).
Another one is given in Proposition 2.2 and Theorem 2.3.

\begin{proposition} (\cite{He3} Prop. 4.1)
The product $(M_1\times M_2,\circ_1\oplus \circ_2,e_1+e_2)$ of two
F-manifolds is an F-manifold. The sum (of the lifts to $M_1\times M_2$)
of two Euler fields $E_i$, $i=1,2,$ on $(M_i,\circ_i,e_i)$ of the
same weight $c\in \C$ is an Euler field of weight $c\in \C$ on 
$M_1\times M_2$.
\end{proposition}

Theorem 2.3 describes the decomposition of a germ of an F-manifold.
In order to 
state it properly we need the following classical and elementary fact:
each tangent space of an F-manifold decomposes as an algebra 
$(T_tM,\circ,e)$ uniquely into a direct sum
\begin{eqnarray}
(T_tM,\circ,e)=\bigoplus_{k=1}^{l(t)}  ((T_tM)_k,\circ,e_k)
\end{eqnarray}
of local subalgebras $(T_tM)_k$ with units $e_k$ and with 
$(T_tM)_j\circ (T_tM)_k=0$ for $j\neq k$ (cf. e.g. \cite{He3} Lemma 1.1).
One can obtain this decomposition as the simultaneous eigenspace decomposition
of the commuting endomorphisms $X\circ: T_tM\to T_tM$ for $X\in T_tM$.

\begin{theorem} (\cite{He3} Theorem 4.2)
Let $(M,\circ,e)$ be an F-manifold and $t\in M$. The decomposition (2.3)
extends to a unique decomposition
\begin{eqnarray}
(M,t)= \prod_{k=1}^{l(t)}(M_k,t)
\end{eqnarray}
of the germ $(M,t)$ into a product of germs of F-manifolds. 
These germs are irreducible germs of F-manifolds as already the algebras
$T_t(M_k) \cong (T_tM)_k$
are irreducible (as they are local algebras).

An Euler field of weight $c$ decomposes accordingly.
\end{theorem}

The proof in \cite{He3} uses (2.1) in a way which justifies calling it 
{\it integrability condition}.

Consider an F-manifold $(M,\circ ,e)$. The function $l:M\to \N$ defined in 
(2.3) is lower semicontinuous (\cite{He3} Proposition 2.3). The {\it caustic}
$\KK := \{ t\ |\ l(t)< \mbox{ generic value }\}$
is empty or a hypersurface (the proof in \cite{He3} Proposition 2.4 for
the case $generic\ value\ =m$  works for any generic value of $l$).

The multiplication on $T_tM$ is {\it semisimple} if $l(t)=m$.
The F-manifold is {\it massive} if the multiplication is generically 
semisimple.

Up to isomorphism there is only one germ of a 1-dimensional F-manifold,
$(\C,\circ, e)$ with $e=\frac{\paa}{\paa u}$ for $u$ a coordinate on 
$(\C,0)$.
The space of Euler fields of weight 0 is $\C\cdot e$, an Euler field of weight
1 is $u\, e$. This germ of an F-manifold is called $A_1$.

By Theorem 2.3, any germ of a semisimple F-manifold is a product 
$A_1^m$, that means, there are local coordinates $u_1,...,u_m$ with
$e_i=\frac{\paa}{\paa u_i}$ and $e_i\circ e_j =\delta_{ij}e_i$.
They are unique up to renumbering and shift and are called canonical 
coordinates, following Dubrovin. The vector fields $e_i$ are called
{\it idempotent}.
Also by Theorem 2.3, then each Euler field of weight 1 takes the form
$\sum_{i=1}^m (u_i+r_i)e_i$ for some $r_i\in \C$.

\begin{example}
Fix $m\geq 1$ and $n\geq 2$. The manifold $M=\C^m$ with coordinate fields
$\delta_i = \frac{\paa }{\paa t_i}$ and multiplication defined by
\begin{eqnarray}
\delta_1\circ \delta_2 &=&\delta_2 \ ,\\
\delta_2\circ \delta_2 &=&t_2^{n-2}\delta_1 \ ,\\
\delta_i\circ \delta_j &=& \delta_{ij}\delta_i \qquad \mbox{ if }
 (i,j)\notin \{ (1,2),(2,1),(2,2)\}
\end{eqnarray}
is a massive F-manifold.
The submanifold $\C^2\times\{ 0\}$ is an F-manifold with the name 
$I_2(n)$, with $I_2(2)=A_1^2$, $I_2(3)=A_2$, $I_2(4)=B_2$, $I_2(5)=H_2$,
and $I_2(6)=G_2$.

$(M,\circ,e)$ decomposes globally into a product 
$\C^2\times \C\times ... \times \C$ of F-manifolds of the type
$I_2(n)A_1^{m-2}$. The unit fields for the components
are $\delta_1, \delta_3, ... , \delta_m$, the global unit field is 
$e=\delta_1+\delta_3+...+\delta_m$, the caustic is 
$\KK=\{ t\ |\ t_2=0\}$. The idempotent vector fields in a simply connected
subset of $M-\KK$ are
\begin{eqnarray}
e_{1/2}&=& \frac{1}{2} \delta_1 \pm 
         \frac{1}{2} t_2^{-\frac{n-2}{2}} \delta_2\ ,\\
e_i &=& \delta_i \qquad \mbox{ for } i\geq 3\ ,
\end{eqnarray}
canonical coordinates there are
\begin{eqnarray}
u_{1/2} &=& t_1 \pm {2\over n} t_2^{n\over 2}\ , \\
u_i &=& t_i \qquad \mbox{ for } i\geq 3\ .
\end{eqnarray}
An Euler field of weight 1 is
\begin{eqnarray}
E=t_1\delta_1 + \frac{2}{m} t_2\delta_2 + \sum_{i\geq 3} t_i\delta_i\ .
\end{eqnarray}
The space of global Euler fields of weight 0 is
$\sum_{i\neq 2} \C\cdot \delta_i$.
\end{example}

The classification of 3-dimensional irreducible germs of massive F-manifolds 
is already vast (\cite{He3} chapter 20). But the classification of 
2-dimensional irreducible germs of massive F-manifolds is nice
(\cite{He3} Theorem 12.1): they are precisely the germs at 0 for 
$m=2$ and $n\geq 3$ in Example 2.4 with the names $I_2(n)$.

\section{Frobenius manifolds}
\setcounter{equation}{0}

Frobenius manifolds were defined first by Dubrovin \cite{Du1}.
They turn up now at many places, see \cite{Du2} and \cite{Man}
(also for more references), especially in quantum cohomology 
and mirror symmetry.
In this paper we will only be concerned with the Frobenius manifolds
in singularity theory (chapter 4).

\begin{definition}
A Frobenius manifold $(M,\circ,e,E,g)$ is a manifold $M$ with a  
multiplication $\circ$ on the tangent bundle, a global unit field $e$,
another global vector field $E$, which is called {\it Euler field}, 
and a metric $g$, subject to the following conditions:
\begin{list}{}{}
\item[1)] the metric is {\it multiplication invariant}, 
         $g(X\circ Y,Z)=g(X,Y\circ Z)$,
\item[2)] (potentiality) the $(3,1)$-tensor $\nabla \circ $ is symmetric
         (here $\nabla$ is the Levi-Civita connection of the metric),
\item[3)] the metric $g$ is flat,
\item[4)] the unit field $e$ is flat, $\nnn e=0$,
\item[5)] the Euler field satisfies $\Lie_E(\circ )=1\cdot \circ$ and 
          $\Lie_E(g)=D\cdot g$ for some $D\in \C$.
\end{list}
\end{definition}

\begin{remarks}
a) Condition 2) implies (2.1) (\cite{He3} Theorem 5.2). Therefore
a Frobenius manifold is an F-manifold.

\smallskip
b) The $(3,0)$-tensor $A$ with $A(X,Y,Z):=g(X\circ Y,Z)$ is symmetric by 1).
Then 2) is equivalent to the symmetry of the $(4,0)$-tensor $\nnn A$.
If $X, Y, Z, W$ are local flat fields then 2) is equivalent to the 
symmetry in $X, Y, Z, W$ of 
$X\, A(Y,Z,W)$. This is equivalent to the
existence of a local {\it potential} $\Phi\in \OO_{M}$ with 
$A(X,Y,Z)=XYZ(\Phi)$ for flat fields $X,Y,Z$.

\smallskip
c) $\Lie_E(g)=D\cdot g$ means that $E$ is a sum of an infinitesimal
dilation, rotation and shift. 
Therefore $\nnn E$ maps a flat field $X$ to a flat field
$\nnn_X E = [X,E] =-\Lie_E X$, i.e. it is a flat
$(1,1)$-tensor, $\nnn (\nnn E)=0$.
Its eigenvalues are called $d_1,...,d_m$.
Now $\Lie_E(\circ )=1\cdot \circ$ implies $\nnn_e E= [e,E]=e$, 
and $\nnn E - \frac{D}{2}\id$ is an infinitesimal isometry because of
$\Lie_E(g)=D\cdot g$.
One can order the eigenvalues such that $d_1=1$ and $d_i+d_{m+1-i}=D$.
\end{remarks}

\section{Hypersurface singularities}
\setcounter{equation}{0}

Let $f:(\C^{n+1},0)\to (\C,0)$ be a holomorphic function germ with an 
isolated singularity at 0 and with Milnor number $\mu$.

An unfolding of $f$ is a holomorphic function germ 
$F:(\C^{n+1}\times \C^m,0)\to (\C,0)$ with
$F|(\C^{n+1}\times \{0\},0) =f$. The coordinates on 
$(\C^{n+1}\times \C^m,0)$ are called $(x_0,...,x_{n},t_1,...,t_m)$.

The germ $(C,0)\subset (\C^{n+1}\times \C^m,0)$ of the critical space
is defined by the ideal $J_F = (\frac{\paa F}{\paa x_0}, ...,
\frac{\paa F}{\paa x_n})$. 
The projection $pr_{C}:(C,0)\to (\C^m,0)$ is finite and flat of degree
$\mu$. The map 
\begin{eqnarray}
\aaa : \TT_{\C^m,0} &\longrightarrow&  \OO_{C,0}\\
\frac{\paa }{ \paa t_i} &\mapsto& \frac{\paa F}{ \paa t_i}|_{(C,0)}
\end{eqnarray}
is the Kodaira-Spencer map.

The unfolding is semiuniversal iff $\aaa$ is an isomorphism.
Consider a semiuniversal unfolding.
Then $m=\mu$, and we set $(M,0):=(\C^\mu,0)$. The map $\aaa$ induces a 
multiplication $\circ$
on $\TT_{M,0}$ by $\aaa(X\circ Y)=\aaa (X)\cdot \aaa (Y)$
with unit field $e=\aaa^{-1}([1])$ and a vector field
$E:=\aaa^{-1}([F])$.

This multiplication and the field $E$ were considered first by K. Saito
\cite{SK1}\cite{SK3}. Because the Kodaira-Spencer map $\aaa$ behaves
well under morphisms of unfoldings, the tuple $((M,0),\circ, e, E)$
is essentially independent of the choices of the semiuniversal unfolding.
This and the following fact are discussed in \cite{He3} (chapter 16).

\begin{theorem}
The base space $(M,0)$ of a semiuniversal unfolding $F$ of an isolated
hypersurface singularity $f$ is a germ $((M,0),\circ,e,E)$ of a 
massive F-manifold with Euler field $E=\aaa^{-1}([F])$. 
For each $t\in M$ there is a canonical isomorphism
\begin{eqnarray}
(T_tM,\circ, E|_t)\cong (\bigoplus_{x\in  Sing (F_t)}
\mbox{ Jacobi algebra of }(F_t,x),\mbox{mult.}, [F_t])\ .
\end{eqnarray}
At generic points of the caustic the germ of the F-manifold is of the 
type $A_2\, A_1^{\mu-2}$.

The base spaces of two semiuniversal unfoldings are canonically isomorphic
as germs of F-manifolds with Euler fields.
\end{theorem}

By work of K. Saito \cite{SK1}\cite{SK3} and M. Saito \cite{SM} one can even
construct a metric $g$ on $M$ such that $((M,0),\circ,e,E,g)$ is  
the germ of a Frobenius manifold.
The construction uses the Gau{\ss}-Manin connections for $f$ and $F$,
K. Saito's higher residue pairings, a polarized mixed Hodge structure,
and results of Malgrange on deformations of microdifferential systems.

A more elementary and much broader version of the construction,
which does not use Malgrange's results, is given in \cite{He4}.
Here we restrict ourselves to a formulation of the result. This uses
the {\it spectral numbers} of $f$ \cite{St1}\cite{AGV}.

Let $f:X\to \Delta$ be a representative of the germ $f$ as usual, with
$\Delta = B_\delta \subset \C$ and
$X=f^{-1}(\Delta)\cap B_\varepsilon^{n+1}\subset \C^{n+1}$
(with $1\gg \varepsilon \gg \delta >0$).
The cohomology bundle 
$H^n = \bigcup_{t\in \Delta^{*}} H^n(f^{-1}(t),\C)$ is flat. 
Denote by $\Delta^\infty \to \Delta^{*}$ a universal covering.
A {\it global flat multivalued section} in $H^n$ is a map
$\Delta^\infty\to H^n$ with the obvious properties.
The $\mu$-dimensional space of the global flat multivalued sections in $H^n$ 
is denoted by $H^\infty$. Steenbrink's Hodge filtration
$F^\bullet$ on $H^\infty$ \cite{St1} together with topological data yields a 
polarized mixed Hodge structure on it (see \cite{He1}\cite{He4}
for definitions and a discussion of this).

The {\it spectral numbers} $\alpha_1,...,\alpha_\mu$ of $f$ are $\mu$ rational
numbers with
\begin{eqnarray}
\sharp (i\ |\ \alpha_i=\alpha)
=\dim \Gr_F^{[n-\alpha]}H^\infty_{e^{-2\pi i\alpha}}\ .
\end{eqnarray}
Here $H^\infty_{e^{-2\pi i\alpha}}$ is the generalized eigenspace of the
monodromy on $H^\infty$ with eigenvalue $e^{-2\pi i\alpha}$.
So $e^{-2\pi i\alpha_1},...,e^{-2\pi i\alpha_\mu}$ are the eigenvalues
of the monodromy. The spectral numbers satisfy 
$-1<\alpha_1\leq ... \leq \alpha_\mu <n$ and 
$\alpha_i+\alpha_{\mu+1-i}=n-1$.

Essential for understanding them and for the whole construction of the 
me\-tric $g$ is Varchenko's way to construct a mixed Hodge structure on 
$H^\infty$ with the Gau{\ss}-Manin connection of $f$
(\cite{Va}\cite{AGV}, cf. also \cite{SchSt}\cite{SM}\cite{He1}\cite{He4}).

It turns out that a metric $g$ such that $((M,0),\circ,e,E,g)$ is a Frobenius
manifold is in general not unique. By work of K. Saito and M. Saito 
each choice of a filtration on $H^\infty$
which is {\it opposite} to $F^\bullet$ (see \cite{SM}\cite{He4} for the
definition) yields a metric which gives a Frobenius manifold structure as 
in Theorem 4.2. A more precise statement and a detailed
proof of it is given in \cite{He4} (chapter 6).

\begin{theorem}
One can choose a metric $g$ on the base space $(M,0)$ of a semiuniversal
unfolding of an isolated hypersurface singularity $f$ such that 
$((M,0),\circ,e,E,g)$ is a germ of a Frobenius manifold and 
$\nnn E$ is semisimple with eigenvalues $d_i=1+\alpha_1-\alpha_i$
and with $D=2-(\alpha_\mu-\alpha_1)$.
\end{theorem}

In fact, often one can also find metrics giving Frobenius manifold structures
with $\{d_1,...,d_\mu\}\neq \{1+\alpha_1-\alpha_i\ |\ i=1,...,\mu\}$
(\cite{SM} 4.4, \cite{He4} Remarks 6.7).

\section{Socle field}
\setcounter{equation}{0}

A Frobenius manifold has another distinguished vector field 
besides the unit field and the Euler field. It will be discussed
in this section. We call it the {\it socle field}.
It is used implicitly in Dubrovin's papers and in \cite{Gi}.

Let $(M,\circ, e, g)$ be a manifold with a multiplication $\circ$ on 
the tangent bundle, with a unit field, and with a multiplication
invariant metric $g$. We do not need flatness and potentiality
and an Euler field in the moment.

Each tangent space $T_tM$ is a Frobenius algebra. This means (more or less
by definition) that the splitting (2.3) 
is now a splitting into a direct sum of Gorenstein rings 
(cf. e.g. \cite{He3} Lemma 1.2)
\begin{eqnarray}
T_tM = \bigoplus_{k=1}^{l(t)} (T_tM)_k 
\end{eqnarray}
They have maximal ideals $\mmm_{t,k}\subset (T_tM)_k$ and units 
$e_k$ such that $e=\sum e_k$. They satisfy
\begin{eqnarray}
{(T_tM)_j}\circ (T_tM)_k =\{ 0\} \qquad \mbox{ for } j\neq k\  ,
\end{eqnarray}
and thus
\begin{eqnarray}
g((T_tM)_j, (T_tM)_k)=\{ 0\} \qquad \mbox{ for } j\neq k\ .
\end{eqnarray}
The socle $\Ann_{(T_tM)_k} (\mmm_{t,k})$ is 1-dimensional and has a 
unique generator $H_{t,k}$ which is normalized such that
\begin{eqnarray}
g(e_k,H_{t,k})=\dim (T_tM)_k\ .
\end{eqnarray}
The following lemma shows that the vectors $\sum_k H_{t,k}$ glue to a 
holomorphic vector field, the {\it socle field} of 
$(M,\circ,e,g)$.

\begin{lemma}
For any dual bases $X_1,...,X_m$ and $\www X_1,...,\www X_m$ of $T_tM$,
that means, $g(X_i,\www X_j)=\delta_{ij}$, one has
\begin{eqnarray}
\sum_{k=1}^{l(t)} H_{t,k} = \sum_{i=1}^m X_i\circ \www X_i\ .
\end{eqnarray}
\end{lemma}

{\bf Proof:}
One sees easily that the sum $\sum X_i\circ \www X_i$ is independent of
the choice of the basis $X_1,...,X_m$. One can suppose that $l(t)=1$ and 
that $X_1,...,X_m$ are chosen such that they yield a splitting of the 
filtration $T_tM\supset \mmm_{t,1}\supset \mmm_{t,1}^2\supset ... $.
Then $g(e,X_i\circ \www X_i )=1$ and 
\begin{eqnarray}
g(\mmm_{t,1},X_i\circ \www X_i )= g(X_i\circ \mmm_{t,1},\www X_i)=0\ .
\nonumber
\end{eqnarray}
Thus $X_i\circ \www X_i = {1\over m} H_{t,1}$. \hfill $\qed$

\bigskip
It will be useful to fix the multiplication and vary the metric.

\begin{lemma}
Let $(M,\circ,e,g)$ be a manifold with multiplication $\circ$ on the 
tangent bundle, unit field $e$ and multiplication invariant metric $g$.
For each multiplication invariant metric $\www g$ there exists a unique
vector field $Z$ such that the multiplication with it is invertible
everywhere and for all vector fields $X,Y$
\begin{eqnarray}
\www g (X,Y)=g(Z\circ X,Y)\ .
\end{eqnarray}
The socle fields $H$ and $\www H$ of $g$ and $\www g$ satisfy
\begin{eqnarray}
H=Z\circ \www H\ .
\end{eqnarray}
\end{lemma}

{\bf Proof:}
The situation for one Frobenius algebra is described for example in 
\cite{He3} (Lemma 1.2). It yields (5.6) immediately. (5.7) follows 
from the comparison of (5.4) and (5.6).
\hfill $\qed$

\bigskip
Denote by
\begin{eqnarray}
H_{op}: \tm \to \tm ,\ X\mapsto H\circ X
\end{eqnarray}
the multiplication with the socle field $H$ of $(M,\circ, e,g)$ as above.
The socle field is especially interesting if the multiplication is 
generically semisimple, that means, generically $l(t)=m$. Then
the caustic $\KK = \{ t\in M\ |\ l(t)<m\}$ is the set where the 
multiplication is not semisimple. It is the hypersurface
\begin{eqnarray}
\KK = \det (H_{op})^{-1}(0)\ .
\end{eqnarray}
In an open subset of $M-\KK$ with basis $e_1,...,e_m$ of idempotent
vector fields the socle field is 
\begin{eqnarray}
H=\sum_{i=1}^m {1\over g(e_i,e_i)}e_i\ .
\end{eqnarray}
It determines the metric $g$ everywhere because (5.10) determines the 
metric at semisimple points.

\begin{theorem}
Let $(M,\circ,e,g)$ be a massive F-manifold with multiplication invariant
metric $g$. Suppose that at generic points of the caustic
the germ  of the F-manifold is of the type $I_2(n)A_1^{m-2}$.

Then the function $\det (H_{op})$ vanishes with multiplicity $n-2$ along 
the caustic.
\end{theorem}

{\bf Proof:}
It is sufficient to consider Example 2.4.
A multiplication invariant metric $g$ is uniquely determined by the 
1-form $\varepsilon = g(e,.)$. Because of (5.7) it is sufficient to
prove the claim for one metric. We choose the metric with 1-form
\begin{eqnarray}
\varepsilon (\delta_i )= 1-\delta_{i1}\ .
\end{eqnarray}
The bases $\delta_1,\delta_2,\delta_3,...,\delta_m$ and
$\delta_2,\delta_1,\delta_3,...,\delta_m$ are dual with respect to this
metric. Its socle field is by Lemma 5.1
\begin{eqnarray}
H=2\delta_2+\delta_3 + ... + \delta_m
\end{eqnarray}
and satisfies $\det (H_{op}) = -4t_2^{n-2}$.
\hfill $\qed$

\section{G-function of a massive Frobenius manifold}
\setcounter{equation}{0}

Associated to any simply connected semisimple Frobenius manifold 
is a fascinating and quite mysterious function. Dubrovin and Zhang
\cite{DZ1}\cite{DZ2} called it the G-function and proved the deepest
results on it. But Givental \cite{Gi} studied it, too, and it 
originates in much older work. It takes the form
\begin{eqnarray}
G(t)= \log \tau_I - {1\over 24} \log J
\end{eqnarray}
and is determined only up to addition of a constant. First we
explain the simpler part, $\log J$. Let $(M,\circ,e,E,g)$ be a 
semisimple Frobenius manifold with canonical coordinates 
$u_1,...,u_m$ and flat coordinates $\www t_1,...,\www t_m$. 
Then
\begin{eqnarray}
J= \det ( {\paa \www t_i \over \paa u_j }) \cdot \mbox{constant}
\end{eqnarray}
is the base change matrix between flat and idempotent vector fields.

One can rewrite it with the socle field.
Denote $\eta_i:=g(e_i,e_i)$ and consider the basis $v_1,...,v_m$ of vector
fields with 
\begin{eqnarray}
v_i={1\over \sqrt{\eta_i}}e_i 
\end{eqnarray}
(for some choice of the square roots).
The matrix $\det (g(v_i,v_j))=1$ is constant as is the 
corresponding matrix for the flat
vector fields. Therefore
\begin{eqnarray}
\mbox{constant} \cdot J = \prod_{i=1}^m \sqrt{\eta_i} = 
\det (H_{op})^{-{1\over 2}} \ .
\end{eqnarray}
Here $H=\sum v_i\circ v_i$ is the socle field.

\bigskip
One of the origins of the first part $\log \tau_I$ are isomonodromic 
deformations. The second structure connections and the first structure
connections of the semisimple Frobenius manifold are 
isomonodromic deformations over $\P^1\times M$ of restrictions to
a slice $\P^1\times\{ t\}$. The function $\tau_I$ is their 
{\it $\tau$-function} in the sense of 
\cite{JMMS}\cite{JMU}\cite{JM}\cite{Mal}.
See \cite{Sab} for other general references on this.

The situation for Frobenius manifolds is discussed and put into a 
Hamilto\-nian framework in \cite{Du2} (Lecture 3), \cite{Man} (II \S 2),
and in \cite{Hi}.
The coefficients $H_i$ of the 1-form $d\log \tau_I = \sum H_idu_i$
are certain Hamiltonians and motivate the definition of this 1-form.
Hitchin \cite{Hi} compares the realizations of this for the first
and the second structure connections.

\bigskip
Another origin of the whole G-function comes from quantum cohomology. 
Getzler \cite{Ge} studied the relations between cycles in the moduli
space ${\oooo \MM}_{1,4}$ and derived from it recursion relations for 
genus one Gromov-Witten invariants of projective manifolds and
differential equations for the genus one Gromov-Witten potential.

Dubrovin and Zhang \cite{DZ1} (chapter 6) investigated these differential
equations for any semisimple Frobenius manifold and found that they have
always one unique solution (up to addition of a constant), the G-function.
Therefore in the case of a semisimple Frobenius manifold coming from
quantum cohomology, the G-function is the genus one Gromov-Witten potential.
(Still I find these differential equations mysterious.)

They also proved part of the conjectures in \cite{Gi}
concerning $G(t)$. Finally, they found that the potential of the 
Frobenius manifold (for genus zero) and the G-function (for genus one)
are the basements of full free energies in genus zero and one
and give rise to Virasoro constraints \cite{DZ2}.
Exploiting this  for singularities will be a big task for the future.

\bigskip
In  chapter 7 we need only the definition of 
$\log \tau_I$ and the behaviour of $G(t)$ with respect to the Euler field
and the caustic in a massive Frobenius manifold. 
We have to resume some known formulas related to the canonical 
coordinates of a semisimple Frobenius manifold
(\cite{Du2}, \cite{Man}, also \cite{Gi}). 

The 1-form $\varepsilon = g(e,.)$ is closed and can be written as 
$\varepsilon=d\eta$. One defines
\begin{eqnarray}
\eta_i &:=& e_i\eta = g(e_i,e)=g(e_i,e_i)\ , \\
\eta_{ij} &:=&  e_ie_j \eta = e_i\eta_j=e_j\eta_i\ ,\\
\gamma_{ij} &:=& {1\over 2} 
        {\eta_{ij} \over \sqrt{\eta_i}\sqrt{\eta_j}}\ ,\\
V_{ij} &:=& -(u_i-u_j) \gamma_{ij}\ , \\
d\log \tau_I &:=& {1\over 2} \sum_{i\neq j} {V_{ij}^2 \over u_i-u_j}du_i 
    = {1\over 2} \sum_{i\neq j} (u_i-u_j) \gamma_{ij}^2 du_i \\
  &=& {1\over 8} \sum_{i\neq j} (u_i-u_j) 
        {\eta_{ij}^2 \over \eta_i\eta_j} du_i \ .
\end{eqnarray}

\begin{theorem}
Let $(M,\circ, e, E,g)$ be a semisimple Frobenius manifold with 
global canonical coordinates $u_1,...,u_m$. 

a) The rotation coefficients $\gamma_{ij}$ (for $i\neq j$) satisfy the 
Darboux-Egoroff equations
\begin{eqnarray}
e_k \gamma_{ij} &=& \gamma_{ik}\gamma_{kj} 
\qquad \mbox{ for } k\neq i \neq j \neq k\ ,\\
e \gamma_{ij} &=& 0 \qquad \qquad \qquad \mbox{for } i\neq j \ .
\end{eqnarray}

b) The connection matrix of the flat connection for the basis $v_1,...,v_m$
from (6.3) is the matrix $\Gamma := (\gamma_{ij}d(u_i-u_j))$. 
The Darboux-Egoroff equations are equivalent to the flatness condition
$d\, \Gamma + \Gamma \land \Gamma = 0$.
 
c) The 1-form $d\log \tau_I$ is closed and comes from a function 
$\log \tau_I$.

d) $E(\eta_i)=(D-2)\eta_i$ and $E(\gamma_{ij})=-\gamma_{ij}$.

e) If the canonical coordinates are chosen such that $E=\sum u_ie_i$ 
then the matrix $-(V_{ij})$ is the matrix of the endomorphism
$\Nu = \nabla E - {D\over 2} \id$ on $\tm $ with respect to the basis
$v_1,...,v_m$. 
\end{theorem}

{\bf Proof:}
a)+b) See \cite{Du2} (pp. 200--201) or \cite {Man} (I \S 3). \\
c) This can be checked easily with the Darboux-Egoroff equations. \\
d) It follows from $\Lie_E (g) =D\cdot g$ and from $[e_i,E]=e_i$. \\
e) This is implicit in \cite{Du2} (pp. 200--201). One can check it with
a)+b)+d).
\hfill $\qed$

\bigskip
The endomorphism $\Nu$ is skewsymmetric with respect to $g$ and flat
with eigenvalues $d_i-{D\over 2}$; the numbers $d_i$ can be ordered
such that $d_1=1$ and $d_i+d_{m+1-i}=D$ (cf. Remark 3.2 c)).

\begin{corollary}
(\cite{DZ1} Theorem 3) Suppose that $E= \sum u_ie_i$. Then
\begin{eqnarray}
E \log \tau_I &=& -{1\over 4} \sum_{i=1}^m (d_i -{D\over 2})^2\ ,\\
E \ G(t) &=& -{1\over 4} \sum_{i=1}^m (d_i-{D\over 2})^2 + 
         {m(2-D)\over 48} =: \gamma \ .
\end{eqnarray}
\end{corollary}

{\bf Proof:}
\begin{eqnarray}
E\log \tau_I &=& {1\over 2} \sum_{i\neq j} {u_iV_{ij}^2\over u_i-u_j} 
  =  {1\over 2} \sum_{i<j} V_{ij}^2 \\
&=& -{1\over 4} \sum_{ij} V_{ij}V_{ji} = -{1\over 4} \trace (V^2) 
\nonumber \\
&=& -{1\over 4} \sum_{i=1}^m (d_i - {D\over 2})^2\ .\nonumber
\end{eqnarray}
(6.4) shows $E(J) = m {D-2\over 2} J$. Now (6.14) follows from the 
definition of the G-function.
\hfill $\qed$

\bigskip
If $M$ is a massive Frobenius manifold with caustic $\KK$, one may ask
which  kind of poles the 1-form $d\log \tau_I$ has along $\KK$ and when
the G-function extends over $\KK$.

All the F-manifolds $I_2(n)$ (cf. Example 2.4) are in a natural way
(up to the choice of a scalar) equipped with a metric $g$, such that they get
Frobenius manifolds (\cite{Du2} Lecture 4, cf. \cite{He3} chapter 19).
These Frobenius manifolds are also denoted $I_2(n)$.

In \cite{DZ1} (chapter 6) the G-function is calculated for them
with coordinates $(t_1,t_2)$ on $M=\C^2$ 
and $e={\paa \over \paa t_1}$. It turns out to be
\begin{eqnarray}
G(t)= -{1\over 24} {(2-n)(3-n)\over n} \log t_2\ .
\end{eqnarray}
Especially, for the case $I_2(3)=A_2$ the G-function is $G(t)=0$.
This was checked independently in \cite{Gi}. Givental concluded
that in the case of singularities the G-function of the base space
of a semiuniversal unfolding with some Frobenius manifold structure
extends holomorphically over the caustic. This is a good guess, but it
does not follow from the case $A_2$, because a Frobenius manifold 
structure on a germ of an F-manifold of type $A_2A_1^{m-2}$ for $m\geq 3$
is never the product of the Frobenius manifolds $A_2$ and $A_1^{m-2}$
(the numbers $d_1,...,d_m$ would not be symmetric).
Anyway, it is true, as the following result shows.

\begin{theorem}
Let $(M,\circ,e,E,g)$ be a simply connected massive Frobenius manifold.
Suppose that at generic points of the caustic $\KK$ the germ of the
underlying F-manifold is of type $I_2(n)A_1^{m-2}$ for one 
fixed number $n\geq 3$.

a) The form $d\, \log \tau_I$ has a logarithmic pole along $\KK$ with
residue $-{(n-2)^2\over 16 n}$ along $\KK_{reg}$.

b) The G-function extends holomorphically over $\KK$ iff $n=3$.
\end{theorem}

{\bf Proof:}
Theorem 5.3 and (6.4) say that the form $-{1\over 24} d\log J$ has a 
logarithmic pole along $\KK$ with residue ${n-2\over 48}$ along 
$\KK_{reg}$. This equals ${(n-2)^2\over 16n}$ iff $n=3$. 
So b) follows from a).

It is sufficient to show a) for the F-manifold in Example 2.4, 
equipped with some metric which makes a Frobenius manifold
out of it (we do not need an Euler field here).
Unfortunately we do not have an identity for $d\log \tau_I$ as (5.7) 
for the socle field which would allow to consider only a 
most convenient metric.

We use (2.5) - (2.11) and (6.5) - (6.10) and consider a neighborhood 
of $0\in \C^m=M$. Denote for $j\geq 3$
\begin{eqnarray}
T_{1j} &:=& (u_1-u_j){\eta_{j1}^2 \over \eta_j\eta_1} + 
            (u_2-u_j){\eta_{j2}^2 \over \eta_j\eta_2} \ ,\\
T_{2j} &:=& (u_1-u_j){\eta_{j1}^2 \over \eta_j\eta_1} - 
            (u_2-u_j){\eta_{j2}^2 \over \eta_j\eta_2} \ ,\nonumber \\
T_{12} &:=& (u_1-u_2){\eta_{12}^2 \over \eta_1\eta_2}  
            d(u_1-u_2)\ . \nonumber
\end{eqnarray}
With $\eta_j(0)\neq 0$ for $j\geq 3$, (6.10) and (2.10) one calculates
\begin{eqnarray}
8d\log \tau_I &=& \mbox{ holomorphic 1-form } + T_{12} \nonumber \\
&+& \sum_{j\geq 3} T_{1j}dt_1 + 
\sum_{j\geq 3}T_{2j}t_2^{{n-2\over 2}}dt_2 \\
&-& \sum_{j\geq 3} T_{1j} du_j \ . \nonumber
\end{eqnarray}
From (2.8) one obtains
\begin{eqnarray}
\eta_{1/2} &=& {1\over 2} \delta_1(\eta) 
         \pm {1\over 2} \delta_2(\eta) t_2^{-{n-2\over 2}}\ ,\\
\eta_1\cdot \eta_2 &=& {1\over 4} t_2^{-n+2} 
        (-\delta_2(\eta )^2 + t_2^{n-2}\delta_1(\eta)^2)\ , \\
\eta_{12} &=& {1\over 4} \delta_1 \delta_1 (\eta) 
+ {1\over 4} {n-2\over 2} t_2^{-n+1}\delta_2(\eta) 
- {1\over 4} t_2^{-n+2}\delta_2\delta_2 (\eta)\ .
\end{eqnarray}
The vector $\delta_2|_0$ is a generator of the socle of the subalgebra in 
$T_0M$ which corresponds to $I_2(n)$. Therefore $\delta_2(\eta)(0)\neq 0$.
It is not hard to see with (6.19) and (2.10) that the terms 
$T_{1j}$ and $T_{2j}t_2^{{n-2\over 2}}$ for $j\geq 3$ are holomorphic at 0.
The term $T_{12}$ is
\begin{eqnarray}
T_{12} &=& {4\over n} \cdot t_2^{n\over 2} \cdot 
{\eta_{12}^2\over \eta_1\eta_2}
\cdot d({4\over n}\cdot t_2^{n\over 2}) \nonumber \\
&=& {8\over n} \cdot t_2^{n-1} \cdot {\eta_{12}^2\over \eta_1\eta_2} \cdot 
dt_2 \\
&=& -{(n-2)^2\over 2n} \cdot {dt_2\over t_2} + \mbox{ holomorphic 1-form}\ .
             \nonumber 
\end{eqnarray}
This proves part a).
\hfill $\qed$

\begin{remark}
It might be interesting to look for massive Frobenius manifolds
which meet the case $n=3$ in Theorem 6.3, but where the underlying
F-manifolds are not locally products of those from hypersurface
singularities. In view of \cite{He3} (Theorem 16.6) the {\it analytic spectrum}
$\Specan (\tm, \circ) \subset T^{*}M$ 
of such F-manifolds would have singularities, but only in codimension
$\geq 2$, as the analytic spectrum of $A_2$ is smooth. 

The analytic spectrum is Cohen-Macaulay and even Gorenstein and a 
Lagrange variety (\cite{He3} chapter 6). P. Seidel (Ecole Polytechnique) 
showed me a normal
and Cohen-Macaulay Lagrange surface. But it seems to be unclear 
whether there exist normal and Gorenstein Lagrange varieties
which are not smooth.
\end{remark}

\section{Variance of the spectrum}
\setcounter{equation}{0}

By Theorem 6.3 the germ $(M,0)$ of a Frobenius manifold as in Theorem 4.2
for an isolated hypersurface singularity $f$ has a holomorphic G-function
$G(t)$, unique up to addition of a constant. 
By Corollary 6.2 and Theorem 4.2 this function satisfies
\begin{eqnarray}
E \ G(t) = -{1\over 4}\sum_{i=1}^\mu (\alpha_i -{n-1\over 2})^2 +
{\mu(\alpha_\mu -\alpha_1)\over 48} =:\gamma \ .
\end{eqnarray}
So it has a very peculiar strength: it gives a grip at the squares of the
spectral numbers $\alpha_1,...,\alpha_\mu$ of the singularity.
Because of the symmetry $\alpha_i+\alpha_{\mu+1-i}=n-1$, the spectral
numbers are scattered around their {\it expectation value}
${n-1\over 2}$. One may ask about their {\it variance}
${1\over \mu }\sum_{i=1}^\mu (\alpha_i-{n-1\over 2})^2$.

\begin{conjecture}
The variance of the spectral numbers of an isolated hypersurface singularity
is 
\begin{eqnarray}
{1\over \mu}\sum_{i=1}^\mu (\alpha_i-{n-1\over 2})^2 \leq 
{\alpha_\mu -\alpha_1\over 12}\ ,
\end{eqnarray}
or, equivalently,
\begin{eqnarray}
\gamma \geq 0\ .
\end{eqnarray}
\end{conjecture}

\begin{theorem}
In the case of a quasihomogeneous singularity $f$
\begin{eqnarray}
{1\over \mu}\sum_{i=1}^\mu (\alpha_i-{n-1\over 2})^2 =
{\alpha_\mu -\alpha_1\over 12}\ ,
\end{eqnarray}
and
\begin{eqnarray}
\gamma = 0\ .
\end{eqnarray}
\end{theorem}

{\bf Proof:}
$(\OO/J_f, mult.,[f]) \cong (T_0M,\circ, E|_0)$. Here one has $f\in J_f$
and $E|_0=0$ and therefore $E\ G(t)=0$.
\hfill $\qed$

\begin{lemma}
The number $\gamma$ of the sum $f(x_0,...,x_n)+g(y_0,...,y_m)$
of two singularities $f$ and $g$ satisfies
\begin{eqnarray}
\gamma (f+g) = \mu(f)\cdot \gamma (g) + \mu(g)\cdot \gamma(f)\ .
\end{eqnarray}
\end{lemma}

{\bf Proof:}
Let $\alpha_1,...,\alpha_{\mu(f)}$ and $\beta_1,...,\beta_{\mu(g)}$
denote the spectral numbers of $f$ and $g$. Then the spectrum of 
$f+g$ as an unordered tuple is \cite{AGV}\cite{SchSt}
\begin{eqnarray}
(\alpha_i+\beta_j+1\ |\ i=1,...,\mu(f),\ j=1,...,\mu(g))\ .
\end{eqnarray}
This and the symmetry of the spectra yields (7.6).
\hfill $\qed$

\begin{remarks}
a) The only unimodal or bimodal families of not semiquasihomogeneous
singularities are the cusp singularities $T_{pqr}$ and the 8 bimodal
series. The spectral numbers are given in \cite{AGV}. One finds 
\begin{eqnarray}
\gamma(T_{pqr}) = {1\over 24} (1 -{1\over p}-{1\over q}-{1\over r})\geq 0
\end{eqnarray}
with equality only for the simple elliptic singularities.  In the case of the
8 bimodal families one obtains
\begin{eqnarray}
\gamma = \frac{p}{48\cdot \kappa} \cdot \left(1-\frac{1}{p+\kappa}
\right) \geq 0
\end{eqnarray}
with $\kappa := 9, 7, 6, 6, 5$ for $E_{3,p}$, $Z_{1,p}$, $Q_{2,p}$, $W_{1,p}$,
$S_{1,p}$, respectively, and
\begin{eqnarray}
\gamma = \frac{p}{48\cdot \kappa} \cdot \left(1+\frac{1}{p+2\kappa}
\right) \geq 0
\end{eqnarray}
with $\kappa := 6, 5, \frac{9}{2}$ for $W_{1,p}^\sharp$,
$S_{1,p}^\sharp$, $U_{1,p}$, respectively.

\smallskip
b) Checking $\gamma=0$ for the $A_\mu$-singularities is easy. With Lemma 7.3
one obtains immediately $\gamma=0$ for all Brieskorn-Pham singularities.
But this is far from a general elementary proof of Theorem 7.2 for 
all quasihomogeneous singularities.

\smallskip
c) In \cite{SK2} K. Saito studied the distribution of the spectral numbers
and their characteristic function
\begin{eqnarray}
\chi_f := {1\over \mu} \sum_{i=1}^\mu T^{\alpha_i+1}
\end{eqnarray}
heuristically and formulated several questions about them.
The G-function might help to go on with these problems.

\smallskip
d) In the case of a quasihomogeneous singularity with weights 
$w_0,...,w_n\in (0,{1\over 2}]$ and degree 1 the characteristic function
is 
\begin{eqnarray}
\chi_f = {1\over \mu} \prod_{i=0}^n {T-T^{w_i} \over T^{w_i}-1}\ ,
\end{eqnarray}
as is well known. It follows easily from (1.1) and (1.2).

\smallskip
e) One can speculate that the Conjecture 7.1, if it is true, comes from
a deeper hidden interrelation between the Gau{\ss}-Manin connection and
polarized mixed Hodge structures.

In \cite{He4} (Remark 6.7 b)) an example of M. Saito (\cite{SM} 4.4) 
is sketched which leads for the semiquasihomogeneous singularity
$f=x^6+y^6+x^4y^4$ to Frobenius manifold structures with
$\{d_1,..,d_\mu\}\neq\{1+\alpha_1-\alpha_i\ |\ i=1,...,\mu\}$.
The number $\gamma$ in that case is $\gamma = -{1\over 144} <0$ .

\smallskip
f) In the case of the simple singularities $A_k, D_k, E_6, E_7, E_8,$
the parameters $t_1,...,t_\mu$ of a suitably chosen unfolding are 
weighted homogeneous with positive degrees with respect to the Euler 
field. Therefore $G=0$ in these cases (cf. \cite{Gi}).
\end{remarks}


\begin{thebibliography}{SaK5}
\bibitem[AGV]{AGV} Arnold, V.I., S.M. Gusein-Zade, A.N. Varchenko: 
    \quad Singularities of differentiable maps, volume II. Birkh\"auser, 
    Boston 1988.
\bibitem[Du1]{Du1} Dubrovin, B.:\quad Integrable systems in topological 
    field theory. Nucl. Phys. {\bf B 375} (1992), 627--685.  
\bibitem[Du2]{Du2} Dubrovin, B.:\quad  Geometry of 2D topological field 
    theories. In: Integrable sytems and quantum groups. Montecatini, Terme
    1993 (M Francoviglia, S. Greco, eds.).
    Lecture Notes in Math. 1620, Springer-Verlag 1996, pp. 120--348.
\bibitem[DZ1]{DZ1} Dubrovin, B., Y. Zhang:\quad Bihamiltonian hierarchies in 
    2D topological field theory at one-loop approximation.
    Commun. Math. Phys. {\bf 198} (1998), 311--361. 
\bibitem[DZ2]{DZ2} Dubrovin, B., Y. Zhang:\quad 
    Frobenius manifolds and Virasoro
    constraints. Sel. math., New ser. {\bf 5} (1999), 423--466.  
\bibitem[Ge]{Ge}  Getzler, E.: Intersection theory on 
    ${\bar {\mathcal M}}_{1,4}$ and elliptic Gromov-Witten invariants.
    J. Amer. Math. Soc. {\bf 10} (1997), 973--998.
\bibitem[Gi]{Gi}  Givental, A.B.: Elliptic Gromov-Witten invariants and
    the generalized mirror conjecture. In: Integrable systems and 
    algebraic geometry. Proceedings of the Taniguchi Symposium 1997
    (M.-H. Saito, Y. Shimizu, K. Ueno, eds.). World Scientific, River Edge NJ 
    1998, pp. 107--155. 
\bibitem[He1]{He1} Hertling, C.:\quad  Classifying spaces and moduli spaces
    for polarized mixed Hodge structures and for Brieskorn lattices. 
    Compositio Math. {\bf 116} (1999), 1--37.
\bibitem[HM]{HM}   Hertling, C., Yu. Manin:\quad  Weak Frobenius manifolds.
    Int. Math. Res. Notices {\bf 1999--6}, 277--286. 
\bibitem[He3]{He3} Hertling, C.:\quad  Multiplication on the tangent bundle.
    First part of the habilitation (also math.AG/9910116).
\bibitem[He4]{He4} Hertling, C.:\quad Frobenius manifolds and moduli spaces
    for hypersurface singularities. June 2000. Second part of the
    habilitation. 
\bibitem[Hi]{Hi}   Hitchin, N.J.:\quad  Frobenius manifolds (notes by 
    D. Calderbank). In: Gauge Theory and symplectic geometry, Montreal
    1995, J. Hurtubise and F. Lalonde eds., 
    Kluwer Academic Publishers, Netherlands 1997, 69--112.
\bibitem[JM]{JM} Jimbo, M., T. Miwa:\quad Monodromy preserving deformations
    of linear ordinary differential equations with rational coefficients II.
    Physica {\bf 2D} (1981), 407--448.
\bibitem[JMMS]{JMMS} Jimbo, M., T. Miwa, Y. Mori, M. Sato:\quad Density matrix
    of impenetrable Bose gas and the fifth Painlev\'e transcendent.
    Physica {\bf 1D} (1980), 80--158.
\bibitem[JMU]{JMU} Jimbo, M., T. Miwa, K. Ueno:\quad Monodromy preserving 
    deformations of linear ordinary differential equations with rational
    coefficients I. Physica {\bf 2D} (1981), 306--352.
\bibitem[Mal]{Mal} Malgrange, B.:\quad Sur les d\'eformations 
    isomonodromiques, I, II. In: S\'eminaire de'l ENS, 
    Math\'ematique et Physique, 1979--1982,
    Progress in Mathematics vol. 37, Birkh\"auser, Boston 1983, pp. 401--438.
\bibitem[Man]{Man} Manin, Yu.:\quad  Frobenius manifolds, quantum cohomology,
    and moduli spaces. American Math. Society, Colloquium Publ. v. {\bf 47},
    1999.  
\bibitem[Sab]{Sab} Sabbah, C.: D\'eformations isomonodromiques et 
    vari\'et\'es de Frobenius, une introduction. 
    Centre de Mathematiques, Ecole Polytechnique, U.M.R. 7640 du C.N.R.S.,
    no. 2000-05, 251 pages.
\bibitem[SK1]{SK1} Saito, K.:\quad  Primitive forms for a universal 
    unfolding of 
    a function with an isolated critical point, J. Fac. Sci. Univ. Tokyo, 
    Sect. IA Math. {\bf 28} (1982), 775--792.
\bibitem[SK2]{SK2} Saito, K.:\quad The zeroes of characteristic function
    $\chi_f$ for the exponents of a hypersurface isolated singular point.
    In: Algebraic varieties and analytic varieties, Advanced Studies in 
    Pure Math. 1, North-Holland Publ. Company, 1983, pp. 195--217.
\bibitem[SK3]{SK3} Saito, K.:\quad  Period mapping associated to a 
    primitive form.
    Publ. RIMS, Kyoto Univ. {\bf 19} (1983), 1231--1264.
\bibitem[SM]{SM} Saito, M.:\quad On the structure of Brieskorn lattices.  
    Ann. Inst. Fourier Grenoble {\bf 39} (1989), 27--72.
\bibitem[SchSt]{SchSt} Scherk, J., J.H.M. Steenbrink:\quad On the mixed Hodge 
    structure on the cohomology of the Milnor fibre. 
    Math. Ann. {\bf 271} (1985), 641-665.
\bibitem[St1]{St1} Steenbrink, J.H.M.:\quad Mixed Hodge structure on the
    vanishing cohomology. In: Real and complex singularities, Oslo 1976, 
    P. Holm (ed.).Alphen aan den Rijn: Sijthoff and Noordhoff 1977, 
    525--562.
\bibitem[St2]{St2} Steenbrink, J.H.M.:\quad Intersection form for 
    quasi-homogeneous singularities.
    Comp. Math. {\bf 34} (1977), 211--223.
\bibitem[Va]{Va} Varchenko, A.N.:\quad The asymptotics of holomorphic forms 
    determine a mixed Hodge structure. 
    Sov. Math. Dokl. {\bf 22} (1980), 772--775.
\end{thebibliography}
\end{document}